\documentclass{article}
\usepackage{emsnewsletter}

\usepackage[utf8]{inputenc} 
\usepackage[T1]{fontenc}    

\usepackage{url}            
\usepackage{booktabs}       
\usepackage{amsfonts}       
\usepackage{nicefrac}       
\usepackage{microtype}      
\usepackage[dvipsnames]{xcolor}
\usepackage{amsmath}
\usepackage{amssymb}
\usepackage{mathtools}
\usepackage{amsthm}
\usepackage{makecell}
\usepackage[font=small]{caption}
\usepackage{multirow}
\usepackage{makecell}
\usepackage{algorithm}
\usepackage{algorithmic}
\usepackage{lscape}
\usepackage{wrapfig}
\usepackage{makecell}
\usepackage{tabu, booktabs}

\usepackage{listings}
\usepackage{verbatim}
\usepackage[breaklinks=true]{hyperref}%

\allowdisplaybreaks
\usepackage{colortbl}
\definecolor{bgcolor}{rgb}{0.8,1,1}
\definecolor{bgcolor2}{rgb}{0.8,1,0.8}
\usepackage{threeparttable}

\usepackage[capitalize,noabbrev]{cleveref}

\theoremstyle{plain}

\usepackage[]{todonotes}

\def\R{\mathbb{R}}
\def\W{\mathcal W}
\def\u{\boldsymbol u}
\def\C{\mathcal C}

\def\R{\mathbb R}
\def\lm{\lambda}

\def\p{\boldsymbol p}
\def\q{\boldsymbol q}

\def\Z{\mathbb Z}

\def\e{\varepsilon}
\def\la{\langle}
\def\ra{\rangle}
\def\vp{f}
\def\y{\mathbf{y}}

\def\one{{\mathbf 1}}
\def\null{{\mathbf 0}}

\newcommand{\norm}[1]{\left\| #1 \right\|}

\newcommand{\bm}[1]{\mathbf{#1}}

\newcommand{\argmin}{\operatornamewithlimits{argmin}}
\newcommand{\argmax}{\operatornamewithlimits{argmax}}
\newcommand{\diag}{\operatornamewithlimits{diag}}

\def\<#1,#2>{\langle #1,#2\rangle}

\def\dm#1{{\color{black}#1}} 
\def\pd#1{{\color{black}#1}} 
\author{Nazarii Tupitsa (MIPT, IITP RAS, Moscow, Russia), Pavel Dvurechensky (WIAS, Berlin, Germany),  Darina Dvinskikh (HSE University, Moscow, Russia),
Alexander Gasnikov (MIPT, Moscow, Russia; IITP RAS, Moscow, Russia; Caucasus Mathematical Center, Adyghe State University, Maikop, Russia)}
\title{Numerical Methods for \pd{Large-Scale} Optimal Transport}

\begin{document}
   
\maketitle
\tolerance=9999
\begin{abstract}
\pd{The optimal transport (OT) problem is a classical optimization problem having the form of linear programming. Machine learning applications put forward new computational challenges in its solution. In particular, the OT problem defines a distance between real-world objects such as images, videos, texts, etc., modeled as probability distributions. In this case, the large dimension of the corresponding optimization problem does not allow applying classical methods such as network simplex or interior-point methods.
This challenge was overcome by introducing entropic regularization and using the efficient Sinkhorn's algorithm to solve the regularized problem. A flexible alternative is the accelerated primal-dual gradient method, which can use any strongly-convex regularization. We discuss these algorithms and other related problems such as approximating the Wasserstein barycenter together with efficient algorithms for its solution, including decentralized distributed algorithms.}
\end{abstract}

\section{Monge--Kantorovich Problem}
\pd{
The optimal transport (OT) problem  is a particular case of \textit{linear programming} (LP) problem. Linear programming is a branch of mathemetical programming concerning minimization (or maximization) problems with linear objectives and linear constraints.
Pioneering contributions in this are were made by the Soviet mathematician and economist L.V. Kantorovich \cite{kantorovich1960mathematical},  who was the first to discover in 1939 that a wide class of the most important industrial production problems can be described mathematically and solved numerically in the form of linear program. 
Important contributions to the development of linear programming were made by T. Koopmans, G.B. Danzig, who proposed the simplex method in 1949, and I.I. Dikin, who initiated the theory of interior point methods in 1967.
}

The history of optimal transport begins with the work of French mathematician  G. Monge~\cite{monge1781memoire}, who proposed a complicated theory to describe    an optimal mass transportation in a geometric way. Inspired by the problem of  resource allocation, L.V. Kantorovich \cite{kantorovich1942translocation} introduced a relaxation of the OT problem which  allowed him to formulate it as a linear programming problem, and, as a consequence, to apply linear programming methods to solve it. \pd{The main idea was to relax the requirement of deterministic nature of transportation, i.e., that a mass from a source point can only be relocated to one target point,  and to introduce  a probabilistic transport, meaning that a mass from a source point can be split between several target points.} 

Consider the OT problem between two  discrete probability measures $p\in S_n(1)$ and  $q\in S_n(1)$, where $ S_n(1):=\{s \in \mathbb{R}^n_+: \la s, \mathbf{1}\ra = 1\}$ is the standard probability simplex. All the results can be generalized to the case when the supports of these distributions have different sizes.
In the Kantorovich's formulation the Monge map is changed to a coupling matrix, also called transportation plan, $\pi \in \mathbb{R}^{n\times n}_+$ with the elements $\pi_{ij}$ prescribing the amount of mass moved from the source point $i$ to the target point $j$. Admissible transportation plans form the transportation polytope $U(p,q)$ of all coupling matrices with marginals equal to the source $p$ and target $q$.  Formally,
\[
U(p,q) \triangleq\{ \pi\in \R^{n\times n}_+: \pi \one_{n} =p, \pi^T \one_{n} = q\}.
\]
The next important element is the symmetric transportation cost matrix $C  \in \R^{n\times n}_+$ with the elements $C_{ij}$ giving the cost of transportation of a unit of mass from the source point $i$ to the target point $j$.
All this naturally leads to the Monge--Kantorovich problem of finding a transportation plan $\pi$ that minimizes the total cost of transportation of the distribution $p$ to the distribution $q$: 
\begin{equation}\label{def:optimaltransport}
    \W(p, q) \coloneqq \min_{\pi \in U(p,q)} \la C, \pi \ra,
\end{equation}
where
$\left\langle A, B\right\rangle  =   \sum_{i,j}^{n,n} A_{ij} B_{ij}$ denotes the Frobenius inner product of $A$ and $B$.
In the particular case when $C_{ij}=\rho(x_i,y_j)^r$ for some $r\geqslant 1$, metric $\rho$, and $x_i$, $y_j$, $i,j=1,...,n$ being the support points of $p$ and $q$ respectively, $(\W(p, q))^{1/r}$ defines the $r$-Wasserstein distance.


\pd{Clearly, \eqref{def:optimaltransport} is a linear programming problem that can be solved by interior-point methods or the simplex method.
The best known complexity bound by interior-point methods is $\widetilde O \left(n^{5/2}\right)$  \cite{lee2014path}. However, this method requires as a subroutine a practical implementation of the Laplacian linear system solver, which is not yet available in the literature. The authors of~\cite{pele2009fast} proposed an alternative, implementable interior point method for OT with a complexity bound $\widetilde O\left(n^3\right)$. The same complexity is achieved by the network simplex algorithm~\cite{tarjan1997dynamic,gabow1991faster}. 
The complexity bound $\widetilde O\left(n^{5/2}\right)$ was obtained in~\cite{hopcroft1973n} for the assignment problem or the problem of perfect matchings in bipartite graphs. When $n$ is large such complexity becomes prohibitive and in what follows we mainly rely on the regularization techniques initiated by \cite{cuturi2013sinkhorn}. 
For further details of computational optimal transport, not covered in this article, we recommend~\cite{peyre2019computational}. This book reviews OT with a bias toward numerical methods and its applications in data sciences, and sheds light on the theoretical properties of OT that make it particularly useful for some of these applications.
 }

\section{\pd{Solving the OT Problem in Large Dimensions}}
\pd{The success of OT in applications, especially in Machine Learning, is driven by the idea of entropic regularization \cite{cuturi2013sinkhorn}.  
Solving linear programs via entropic regularization may seem difficult since the entropy function is not self-concordant, which makes it not a good barrier function for interior point methods~\cite{boyd2004convex}. Nevertheless, the entropic regularization is widely applied to various linear and nonlinear
problems with empirical success~\cite{fang1997entropy}.
}

\pd{In this section for the most part, the main focus is on finding an $\e$-additive approximate solution to ~\eqref{def:optimaltransport}. Formally, the goal is to find $\hat \pi \in U(p,q)$ s.t.
\begin{equation}\label{OT_primal}
     \la C, \hat \pi \ra -  \la C, \pi^* \ra \leqslant \e,
\end{equation}
where $\pi^*$ is a solution of~\eqref{def:optimaltransport}. To do that, one can add the regularization by the entropy of the transportation plan and consider the following \textit{entropy-regularized OT problem} 
\begin{align}\label{def:regoptimaltransport}
    \W_\gamma(p, q) \coloneqq \min_{\pi \in  U(p, q)} \left\{g(\pi) =\la \pi, C \ra + \gamma H(\pi)\right\},
\end{align}}
\dm{where $H(\pi) = \la \pi, \ln \pi \ra$  is the negative entropy and $\gamma>0$ is the regularization parameter.  It has recently gained popularity in the literature of image processing and machine learning and has been considered to compute smoothed Wasserstein barycenters of discrete measures \cite{ferradans2014regularized}. As $H(\pi)$ is 1-strongly convex on $S_n(1)$ in the  $\ell_1$-norm, the objective in \eqref{def:regoptimaltransport} is $\gamma$-strongly convex  with respect to $\pi$ in the $\ell_1$-norm on $S_n(1)$, and hence problem \eqref{def:regoptimaltransport} has a unique optimal solution. Moreover, $\W_\gamma (p,q)$ is $\gamma$-strongly convex with respect to $p$ in the $\ell_2$-norm on $S_n(1)$ \cite{bigot2019data}, see (Theorem 3.4).}



Let us consider the primal and dual optimization problems 
\begin{equation}\label{primal_general}
    \min_{\pi \in Q \subseteq E, \,\, \bm{A}\pi=b }  g(\pi),
\end{equation}
\begin{equation}\label{dual_general}
    \max_{\mu \in \Lambda} \left\{   \phi(\mu)= -\la \mu, b \ra  + \min_{\pi \in Q} \left( g(\pi) + \la \bm{A}^T \mu, \pi \ra \right) \right\},
\end{equation}
where $E$ is a finite-dimensional real vector space, $Q$ is a simple closed convex set, $g$ is convex function, $\bm{A}$ is a given linear operator from $E$ to some finite-dimensional real vector space $H$, $\Lambda=H^*$ is the conjugate space. In the particular case of OT problem $\bm{A}$ is s.t. $\bm{A}\pi = \left(\left(\pi \one_{n}\right)^T, \left(\pi^T \one_{n}\right)^T\right)^T$, $b \in H$, $b = (p^T, q^T)^T$.

Since $g$ is convex, $-\phi(\mu)$ is a convex function and, by Danskin's theorem, its subgradient is equal to
\begin{equation}\label{eq:nvp}
    -\nabla \phi(\mu) = b - \bm{A} \pi (\mu),
\end{equation}
where $\pi (\mu)$ is some solution of the convex problem
\begin{equation}\label{eq:inner}
    \min_{\pi \in Q} \left( g(\pi) + \la \bm{A}^T \mu, \pi \ra \right).
\end{equation}

For the particular formulation of the OT problem~\eqref{def:regoptimaltransport} by introducing $Q = \R_+^{n\times n}$ and the Lagrange multipliers $y \in \R^n$ and $z \in \R^n$ for the linear equality constraints in \eqref{def:regoptimaltransport} the dual problem can be written as
\begin{equation}
    \max_{y,z\in\mathbb{R}^n}  \left\{
    \phi(y,z)= \la y,p \ra + \la z,q \ra + \frac{\gamma}{e}
    \sum_{i,j=1}^n e^{\left(- \frac{y_i+z_j+C_{i j}}{\gamma} \right)}\right\}.
\end{equation}
Making the change of variables
\begin{equation}\label{var_change}
    u=-y / \gamma -1/2, \,\, v=-z / \gamma -1/2,
\end{equation}
one arrives at an equivalent, but more standard formulation of the dual problem as a minimization problem
\begin{equation} \label{OT_dual}
    \min_{u,v\in\mathbb{R}^n} \left\{ \vp(u,v)=\gamma\left(\one_n^T B(u,v)\one_n-\la u,p \ra - \la v,q \ra\right)\right\},
\end{equation}
 where $[B(u,v)]_{ij} = \exp \left(u_{i}+v_{j} -\frac{C_{i j}}{\gamma}\right)$. 
Problem~\eqref{OT_dual} is smooth, convex, and unconstrained, and has a natural coupling to the primal variable
 \begin{equation}\label{eq:primal}
     \pi(u,v) = B(u,v).
 \end{equation}

\pd{The state-of-the-art algorithm for solving the regularized OT problem~\eqref{OT_dual} is the \textbf{Sinkhorn's algorithm or RAS  matrix scaling algorithm}~\cite{sinkhorn1974diagonal,knight2008sinkhorn,kalantari2008complexity}.
}As well as Sinkhorn's algprithm most algorithms for the approximating the OT distance are methods for solving the dual problem. Thus, having an output of such algorithm vector $\mu^t$ one can obtain for $g$ from~\eqref{def:regoptimaltransport}
\begin{multline}
    g\left(\pi\left(\mu^t\right)\right) - g(\pi_\gamma^*) =  g\left(\pi\left(\mu^t\right)\right) - \phi(\mu^*) \leqslant g\left(\pi\left(\mu^t\right)\right) - \phi(\mu)
    \\
    = - \la \mu, A\pi(\mu^t) -b\ra = \la \mu, \nabla \phi(\mu) \ra,
\end{multline}
where $\pi_\gamma^*$ is a solution to the minimization problem in~\eqref{def:regoptimaltransport}.



From~\cite{dvurechensky2018computational} one has that
\[
    \la \mu^t, \nabla \phi(\mu^t) \ra \leqslant \frac{\gamma R}{2} E^t\left(\mu^t\right),
\]
where $R = \|C\|_\infty / \gamma - \ln\left(\min\left\{\min_i p_i, \min_j q_j\right\}\right)$ and $E^t\left(\mu^t\right) =  \|\pi\left(\mu^t\right)\one_n - p\|_1 + \|\left(\pi\left(\mu^t\right)\right)^T\one_n - q\|_1$.

Next, consider the complexity of approximating regularized OT
\begin{equation} \label{prob:reg_ot}
        g\left(\pi\left(\mu^t\right)\right) - g\left(\pi_\gamma^*\right) \leqslant \e.
\end{equation}

The authors of~\cite{dvurechensky2018computational} showed that Sinkhorn's algorithm~\cite{sinkhorn1974diagonal}, that iteratively minimizes~\eqref{OT_dual}, after $t>2$ steps ensures  that 
$E^t\left(\mu^t\right) \leqslant \frac{4R}{t-2}$.
The latter leads to $\widetilde O \left({n^2\|C\|^2_\infty}/{\gamma\e}\right)$ arithmetic operation complexity of approximating~\eqref{prob:reg_ot} in the case $\|C\|_\infty / \gamma \gtrapprox - \ln\left(\min\left\{\min_i p_i, \min_j q_j\right\}\right)$. The same result was proved in \cite{lin2019efficient} for Greenkhorn, a greedy version of Sinkhorn, proposed in~\cite{altschuler2017near-linear}.

The authors of~\cite{lin2019efficiency} proposed an accelerated version of Sinkhorn's algorithm, and proved that $E^t\left(\mu^t\right) \leqslant \e'$ in $t \leqslant 1+\left(\frac{16\sqrt{n}R}{\e'}\right)^{2/3}$. The latter leads to $\widetilde O \left( {n^{7/3}\|C\|^{4/3}_\infty}/{\left(\gamma\e\right)^{2/3}}\right)$ arithmetic operation complexity of approximating~\eqref{prob:reg_ot} in the case $\|C\|_\infty / \gamma \gtrapprox - \ln\left(\min\left\{\min_i p_i, \min_j q_j\right\}\right)$.


    

\pd{
The regularized objective in~\eqref{def:regoptimaltransport} uniformly approximates the objective in~\eqref{def:optimaltransport} since $-H(\pi) \in [-2\ln n,0]$ for any $\pi \in U(p,q)$. This also motivates the choice $\gamma = O\left(\frac{\e}{\ln n}\right)$ since then the approximation error is $O\left(\e\right)$. Such choice of $\gamma$ leads to the following complexity results for solving the non-regularized problem \eqref{def:optimaltransport} via the Sinkhorn's algorithm.
}
In~\cite{altschuler2017near-linear}, the authors 
proved the complexity $\widetilde O \left(n^{2}\|C\|^3_\infty /\e^3\right)$ arithmetic operations for \pd{the Sinkhorn's algorithm} and its greedy variant named Greenkhorn to solve the original problem \eqref{def:optimaltransport} in the sense of \eqref{OT_primal}. The authors of~\cite{dvurechensky2018computational} improved the analysis of \pd{the Sinkhorn's algorithm} and proved the complexity $\widetilde O \left(n^{2}\|C\|^2_\infty /\e^2\right)$. 

\pd{ 
Moreover, Sinkhorn's algorithm alternatively balances the matrix $\pi(u,v)$ such that it alternatively has one marginal equal to $p$ or the other marginal equal to $q$. The authors of~\cite{franklin1989scaling} showed the linear convergence of these iterations, i.e., that they form a contraction for 
$e^{u^t - u^*}$ and $e^{v^t - v^*}$ in Hilbert’s projective metric, where $\left(u^*, v^*\right)$ is the minimum in~\eqref{OT_dual}. The result was slightly improved in~\cite{stonyakin2019gradient}.
}

The complexity~$\widetilde O \left(n^{2}\|C\|^2_\infty /\e^2\right)$ was obtained for Greenkhorn in~\cite{lin2019efficient}.  In~\cite{lin2019efficiency} was proposed an accelerated version of the Sinkhorn's algorithm) that has the complexity~$\widetilde O \left(n^{7/3}\|C\|^{4/3}_\infty /\e^{4/3}\right)$. 
The authors of~\cite{altschuler2018approximating} proposed a simple deterministic algorithm that, given two distributions supported on $n$ points in $\R^d$, can compute an
additive $\e$-approximation to the quadratic transportation cost in time $\widetilde O \left(\frac{n}{\e^3}\left(\frac{B\log(n)}{\e}\right)^d\right)$, where $B >0$ is  a universal constant.

In~\cite{cominetti1994asymptotic}, the authors studied the question of how well the solution\pd{, i.e., the transportation plan,} to the regularized problem approximates solutions  of~\eqref{def:optimaltransport} as a function of $\gamma$ and showed that the optimal solution of~\eqref{def:regoptimaltransport} approaches an optimal solution of~\eqref{def:optimaltransport} exponentially fast. In~\cite{weed2018explicit}, the author improved the result and provided explicit exponential bound with easy-to-understand constants. Moreover, the author of~\cite{weed2018explicit} showed that an $\e$-solution to the problem~\eqref{def:optimaltransport} with cost matrix $C \in \Z^{n\times n}_{\geqslant 0}$ can be found by solving the entropy-regularized  problem~\eqref{def:regoptimaltransport} with parameter $\gamma = O\left(\frac{1}{n\ln n / \e}\right)$.

Besides \pd{the Sinkhorn's algorithm, \textbf{accelerated gradient methods} are adapted for solving OT problems.} These methods achieve better theoretical convergence rates compared to Sinkhorn-like methods in some regimes. To the best of our knowledge, the first such method was proposed in~\cite{gasnikov2016efficient}, where the authors proposed non-adaptive Accelerated Gradient Descent (AGD) method \pd{for a more general class of entropy-linear programming problems that includes problem \eqref{def:regoptimaltransport} as a special case.} The approach is based on introducing a tautological constraint set $Q =\{ \pi \in \R_+^{n \times n}: \one_n^T \pi\one_n = 1\}$ in~\eqref{def:regoptimaltransport}, which makes the objective in the corresponding dual problem having Lipchitz continuous gradient. This eventually leads to the  dual problem 
\begin{equation} \label{OT_dual_lip}
    \min_{u,v\in\mathbb{R}^n}  \left\{\varphi(u,v)=\gamma\left(\ln\left(\one_n^T B(u,v)\one_n\right)-\la u,p \ra - \la v,q \ra\right)\right\}
\end{equation}
with the primal-dual coupling (cf. \eqref{eq:primal})
\begin{equation}  \label{eq:pvr}
    \pi(u,v) = B(u,v) / \one_n^T B(u,v)\one_n.
\end{equation}
\pd{The algorithmic idea is to run AGD for solving \eqref{OT_dual_lip} and equip it with some primal updates to guarantee the convergence rate also for the primal problem.}
Importantly, as observed in \cite{dvurechensky2018computational}, this approach is flexible enough to allow working with different strongly-convex regularizers in~\eqref{def:regoptimaltransport} instead of the entropy, which is not possible for the Sinkhorn's algorithm. In particular in some applications the regularization with the squared Euclidean norm is preferable~\cite{blondel2018smooth}. \pd{Moreover, these algorithms allow solving general linearly-constrained strongly-convex optimization problems.}

The line of research \cite{dvurechensky2018computational,lin2019efficient,guminov2021combination,guo2020fast} on primal-dual AGD for solving the pair of problems~\eqref{def:regoptimaltransport} and~\eqref{OT_dual_lip} eventually provides the complexity $\widetilde O \left({n^{5/2}\|C\|_\infty}/{\left(\gamma\e\right)^{1/2}}\right)$ for regularized OT problem~\eqref{def:regoptimaltransport} in the case $\|C\|_\infty / \gamma \gtrapprox - \ln\left(\min\left\{\min_i p_i, \min_j q_j\right\}\right)$ and $\widetilde O \left(n^{5/2} \|C\|_\infty /\e\right)$ arithmetic operations for solving the original problem \eqref{def:optimaltransport} in the sense of \eqref{OT_primal} with the proper choice of $\gamma = O\left(\frac{\e}{\ln n}\right)$. The authors of \cite{guminov2021combination} proposed \pd{an accelerated alternating minimization (AAM) algorithm}, which is a hybrid of the Sinkhorn's algorithm and primal-dual  AGD with the same complexities $\widetilde O \left({n^{5/2}\|C\|_\infty}/{\left(\gamma\e\right)^{1/2}}\right)$ for regularized OT problem~\eqref{def:regoptimaltransport} in the case $\|C\|_\infty / \gamma \gtrapprox - \ln\left(\min\left\{\min_i p_i, \min_j q_j\right\}\right)$ and $\widetilde O \left(n^{5/2} \|C\|_\infty /\e\right)$ for solving the original problem \eqref{def:optimaltransport} in the sense of \eqref{OT_primal} with the proper choice of $\gamma = O\left(\frac{\e}{\ln n}\right)$.



\pd{
Both, the Sinkhorn's algorithm and AGD use as a basic operation the multiplication of the matrix $\pi(u,v)$ by a vector, which may be costly in large dimensions. As a remedy to this,  the authors of~\cite{genevay2016stochastic} formulated the dual problem in the form of expectation using slightly different entropic penalty $\gamma \sum_{i,j}^{n,n} \pi_{ij}\log \frac{\pi_{ij}}{p_i q_j}$, which allowed to apply stochastic gradient descent methods such as SAG ~\cite{schmidt2017minimizing}, and, thereby, propose an efficient algorithm for computing large-scale optimal transport.
}

\pd{The next describes approaches to solving the original problem \eqref{def:optimaltransport} in the sense of \eqref{OT_primal} that do not use an entropy regularization. 
The authors of~\cite{jambulapati2019direct} proposed a new method based on \textbf{dual-extrapolation}~\cite{nesterov2007dual} algorithm and area convexity~\cite{sherman2017area}, and obtained }  the \pd{complexity} $\widetilde O \left(n^2 \|C\|_\infty/\e\right)$ (see also~\cite{blanchet2018towards} for similar rate). Despite the theoretical optimality of the complexity $O \left(n^2 \|C\|_\infty/\e\right)$ (see~\cite{blanchet2018towards}), AGD or AAM outperform in practice the area-convexity approach as it is shown by the numerical experiments of~\cite{guminov2021combination} and~\cite{dvinskikh2021improved}.
\pd{Close to the previous approach is the celebrated  Primal-Dual Hybrid Gradient algorithm \cite{chambolle2011first-order} designed for solving a very general class of saddle-point problems.} In~\cite{chambolle2022accelerated}, this algorithm was applied to OT problems which resulted in the complexity bound $\widetilde O \left(n^{5/2} \|C\|_\infty/\e\right)$ matching the bounds of existing practical state-of-the-art methods.

It is also possible to use a carefully designed \textbf{Newton-type} algorithm to solve the OT problem~\cite{allen2017much,cohen2017matrix}, by making use of a connection to matrix scaling problems. For solving the original problem~\eqref{def:optimaltransport} in the sense of \eqref{OT_primal}, the authors of~\cite{blanchet2018towards} and~\cite{quanrud2018approximating} provided a complexity bound of $\widetilde O \left(n^2 \|C\|_\infty/\e\right)$ for Newton-type algorithms. However, these methods are complicated and efficient implementations are not yet available.
A variety of standard optimization algorithms have also been adapted to the OT setting, e.g. quasi-Newton methods~\cite{cuturi2016smoothed,blondel2018smooth}, yet without the complexity guarantees.

The complexity bounds of all aforementioned methods is presented in Table~\ref{Tab:entr-reg-ot} and Table~\ref{Tab:ot}.

\begin{table}[ht!]
    \centering
    \caption{ Algorithm for OT problem with entropy penalization}
    \begin{threeparttable}
    \begin{tabular}{ll}
        Algorithm   &     Complexity~\eqref{def:regoptimaltransport}$^{(1)}$ 
        \\
         \midrule
        Sinkhorn~\cite{kroshnin2019complexity}  & $\widetilde O \left({n^2\|C\|^2_\infty}/{\gamma\e}\right)$  
        \\
        Greenkhorn~\cite{lin2019efficient} &    $\widetilde O \left({n^2\|C\|^2_\infty}/{\gamma\e}\right)$      
        \\
        Algorithm 7~\cite{lin2019efficiency}   & $\widetilde O \left( {n^{7/3}\|C\|^{4/3}_\infty}/{\left(\gamma\e\right)^{2/3}}\right)$  
        \\
        AGD~\cite{dvurechensky2018computational,lin2019efficient,guo2020fast}  &  $\widetilde O \left({n^{5/2}\|C\|_\infty}/{\left(\gamma\e\right)^{1/2}}\right)$ 
        \\
        AAM~\cite{guminov2021combination}  &  $\widetilde O \left({n^{5/2}\|C\|_\infty}/{\left(\gamma\e\right)^{1/2}}\right)$ 
    \end{tabular}
    \begin{tablenotes}
    \item[]  $^{(1)}$ Complexity bounds are obtained in the case $\|C\|_\infty / \gamma \gtrapprox - \ln\left(\min\left\{\min_i p_i, \min_j q_j\right\}\right)$.
    \end{tablenotes}    
    \end{threeparttable}
    \label{Tab:entr-reg-ot}
\end{table}

For obtaining the bounds for solving the original problem \eqref{def:optimaltransport} in the sense of \eqref{OT_primal} it is sufficient to substitute $\gamma = O\left(\frac{\e}{\ln n}\right)$.

\begin{table}[ht!]
    \centering
    \caption{ Algorithm for OT problem without entropy penalization}
    \resizebox{\columnwidth}{!}{%
    \begin{threeparttable}
    \begin{tabular}{lll}
        Algorithm   & Complexity~\eqref{def:optimaltransport}   & Implementability \\
         \midrule
        Path finding~\cite{lee2014path} & $\widetilde O \left(n^{5/2}\right)$ & $\times$
        \\
        Interior point~\cite{pele2009fast} & $\widetilde O\left(n^3\right)$ & $\surd$
        \\
        Network simplex~\cite{tarjan1997dynamic,gabow1991faster}& $\widetilde O\left(n^3\right)$ & $\surd$
        \\
        Algorithm 1~\cite{hopcroft1973n}$^{(1)}$& $\widetilde O\left(n^5/2\right)$ & $\surd$
        \\
        Dual-extrapolation~\cite{jambulapati2019direct} &   $\widetilde O \left(n^{2}\|C\|_\infty /\e\right)$ & $\surd$
        \\
        Mirror-prox~\cite{jambulapati2019direct} &   $\widetilde O \left(n^{5/2}\|C\|_\infty /\e\right)$ & $\times$
        \\
        Newton-type~\cite{blanchet2018towards,quanrud2018approximating}  & $\widetilde O \left(n^2 \|C\|_\infty/\e\right)$ & $\times$ 
        \\
        PDHG~\cite{chambolle2022accelerated} &  $\widetilde O \left(n^{5/2} \|C\|_\infty/\e\right)$ & $\surd$ 
        \\
        Algorithm 1~\cite{altschuler2018approximating}$^{(2)}$& $\widetilde O \left(\frac{n}{\e^3}\left(\frac{B\log(n)}{\e}\right)^d\right)$ & $\times$ 
        
    \end{tabular}
    \begin{tablenotes}
    \item[] $\times$ indicates that an algorithm requires a subroutine that practical implementation is not yet available in the literature or an algorithm was not implemented by authors.
    \item[] $\surd$ indicates that method can be implemented in practice.
    \item[]  $^{(1)}$ For OT formulated as an assignment problem.
    \item[]  $^{(2)}$ Result for the quadratic transportation cost.
    \end{tablenotes}    
    \end{threeparttable}
    }
    \label{Tab:ot}
\end{table}

In the next two subsections we describe in detail the Sinkhorn's algorithm and accelerated alternating minimization method for the dual problem~\eqref{OT_dual}.

\subsection{Sinkhorn's Algorithm} \label{sec:sink}

\pd{Sinkhorn's algorithm~\cite{cuturi2013sinkhorn}  is the main workhorse for solving large-scale OT problems in applications. The algorithm sequentially minimizes the dual function in~\eqref{OT_dual} with respect to either $u$ or $v$ variable with $v$ or $u$ fixed respectively.
The following Lemma~\ref{Sinkhorn=AM} shows that the dual objective in~\eqref{OT_dual} can be minimized exactly and explicitly over the variable $u$ or $v$. Algorithm \ref{Sinkhorn} presents the pseudocode of the Sinkhorn's algorithm.}

\begin{lemma} \label{Sinkhorn=AM}
    The iterations of Sinkhorn algorithm
    \[u^{t+1} \in \argmin_{u\in\mathbb{R}^n}\vp(u,v^t),\, v^{t+1}\in \argmin_{v\in\mathbb{R}^n}\vp(u^{t+1},v),\] can be written explicitly as
    \[u^{t+1}=u^{t}+\ln p-\ln \left(\pi\left(u^t, v^t\right) \one_n\right), \]
    \[v^{t+1}=v^{t}+\ln q-\ln \left(\pi\left(u^{t+1}, v^t\right)^T \one_n\right), \]
    where $\pi(u,v)$ is given by~\eqref{eq:primal}.
\end{lemma} 

\begin{algorithm}[H]
    \caption{Sinkhorn's algorithm}
    \label{Sinkhorn}
    \begin{algorithmic}
        \REQUIRE $C$, $p, q$, $\gamma > 0$
        \FOR{$t\geqslant 1$}
            \IF{$t$ mod $2 =0$}
            \STATE $u^{t+1}=u^{t}+\ln p-\ln \left(\pi\left(u^{t}, v^{t}\right) \one_n\right)$, $v^{t+1}=v^{t}$
            \ELSE
            \STATE $u^{t+1}=u^{t}$, $v^{t+1}=v^{t}+\ln q-\ln \left(\pi\left(u^{t}, v^{t}\right)^T \one_n\right)$
            \ENDIF
            \STATE $t=t+1$
        \ENDFOR
         \ENSURE $\check{\pi} = \pi(u^t, v^t)$
    \end{algorithmic}
\end{algorithm}

Following~\cite{leonard2013survey,benamou2015iterative} the problem~\eqref{def:regoptimaltransport} can be refactored in a Kullback--Leibler projection form
\begin{align}\label{prob:bary_KL}
    \min\limits_{\pi \in \C_1 \cap  \C_2}  KL\left(\pi | K\right),
\end{align}
where $KL(\pi | \pi') \coloneqq \sum_{i, j = 1}^n \left(\pi_{i j} \ln \left(\frac{\pi_{i j}}{\pi'_{i j}}\right) - \pi_{i j} + \pi'_{i j}\right) = \la \pi, \ln \pi - \ln \pi' \ra + \la \pi' - \pi, \one_n \one_n^T \ra$, $K = \exp\left(-{C}/{\gamma}\right)$ and the affine convex sets $\C_1$ and $\C_2$ with  
\begin{equation}\label{eq:C1C2Def}
    \C_1 = \left\{\pi  : \pi \one_n = p \right\},  \quad
    \C_2 = \left\{\pi  : \pi^T \one_n = q \right\}. 
\end{equation}
The Sinkhorn's algorithm (Algorithm~\ref{Sinkhorn}) 
consists in alternating projections to the sets $\C_1$ and $\C_2$ w.r.t. Kullback--Leibler divergence. 
This algorithm is equivalent to alternating minimization of the dual problem
~\eqref{OT_dual} to the problem~\eqref{def:regoptimaltransport} and can be refactored as Algorithm~\ref{sink-primal}.

\begin{algorithm}[H]
    \caption{KL projection form Sinkhorn}
    \label{sink-primal}
    \begin{algorithmic}[1]
        \REQUIRE Cost matrix $C$, probability measures $p, q$, $\gamma > 0$, starting transport plans $\pi^0 \coloneqq \exp\left(-\frac{C}{\gamma}\right)$
        \FOR{t=0,1,...}
            \IF{$t \bmod 2 = 0$}
                \STATE $\pi^{t + 1} \coloneqq \argmin\limits_{\pi \in \C_1} KL\left(\pi | \pi^t\right)$
            \ELSE
                \STATE $\pi^{t + 1} \coloneqq \argmin\limits_{\pi \in \C_2} KL\left(\pi | \pi^t\right)$
            \ENDIF
            \STATE $t \coloneqq t + 1$
        \ENDFOR
        \ENSURE $\pi^t$    
    \end{algorithmic}
\end{algorithm}

\pd{Algorithm~\ref{Sinkhorn=AM} is the so-called log-domain (equivalent) version of the Sinkhorn's algorithm that has better numerical stability for small $\gamma$. In practice, especially for large $\gamma$, the matrix-scaling implementation in~\cite{cuturi2013sinkhorn} can be more efficient, in particular, due to its GPU-friendly form.}

The iterations of Algorithm~\ref{Sinkhorn} should be implemented using numerically stable computation of the log of the sum of exponents. The example of Python implementation of the Sinkhorn's iteration is given below:

\begin{small}
\begin{lstlisting}[language=Python]
  from scipy.special import logsumexp
  u=logsumexp((-v-C/gamma).T,b=1/q,axis=0)
  v=logsumexp((-u-C.T/gamma).T,b=1/p,axis=0) 
\end{lstlisting}
\end{small}

The next Theorem~\ref{th:pure_sink} from~\cite{dvurechensky2018computational} states the iteration complexity bound for the Sinkhorn’s algorithm.
\begin{theorem} \label{th:pure_sink}
    After $t>2$ steps of Algorithm~\ref{Sinkhorn=AM} it holds that 
    \begin{equation}
        \|\check{\pi}\one_n - p\|_1 + \|\check{\pi}^T\one_n - q\|_1 \leqslant \frac{4R}{t-2},
    \end{equation}
    where $R = \|C\|_\infty / \gamma - \ln\left(\min\left\{\min_i p_i, \min_j q_j\right\}\right)$.
\end{theorem}

\pd{Since the Sinkhorn's algorithm solves the dual problem \eqref{OT_dual},  the approximate solution $\check{\pi}$ is not guaranteed to belong to the feasible set $U(p, q)$.} The authors of~\cite{altschuler2017near-linear} proposed  Algorithm~\ref{Alg:round}, that rounds (projects) an approximate solution $\check{\pi}$ to ensure the projection to lie in the feasible set, i.e., for a given matrix $\check{\pi}$ it outputs $\hat{\pi} \in U(p, q)$ s.t.
\[
\|\check{\pi} - \hat{\pi}\|_1 \leqslant \|\check{\pi}\one_n - p\|_1 + \|\check{\pi}^T\one_n - q\|_1.
\]

\begin{algorithm}[!ht]
    \caption{Rounding algorithm~ \cite{altschuler2017near-linear}}
    \label{Alg:round}
    \begin{algorithmic}[1]
        \REQUIRE{$\check{\pi}$}
        \STATE $P := \diag\left\{\min\left(p_i/p_i(\check{\pi}), 1\right)\right\}$, where $p_i(\check{\pi}) =[\check{\pi}\one_n]_i$
        \STATE $Q := \diag\left\{\min\left(q_i/q_i(\check{\pi}), 1\right)\right\}$, where $q_i(\check{\pi}) =[\check{\pi}^T\one_n]_i$
        \STATE $F := P\check{\pi}Q$
        \STATE $e_p = p - p(F)$, $e_q = q - q(F)$
        \STATE $\hat{\pi} = F + e_p e_{q}^T / \|e_p\|_1$
        \ENSURE $\hat{\pi}$
  \end{algorithmic}
\end{algorithm}

The iterations of the Sinkhorn's algorithm require element-wise logarithm of measures $p$ and $q$, that should be guaranteed to be separated from zero to make it correctly defined. Moreover, Theorem~\ref{th:pure_sink} states that the smaller $\min\left\{\min_i p_i, \min_j q_j\right\}$ the slower is the convergence rate. So, the authors of~\cite{dvurechensky2018computational} proposed 
to replace the input measures $p$ and $q$ with $\check{p}$ and $\check{q}$ s.t. $\|p - \check{p}\|_1 \leqslant \e'/4$, $\|q - \check{q}\|_1 \leqslant \e'/4$  and $\min_i \check{p}_i \geqslant \e'/8n$, $\min_i \check{q}_i \geqslant \e'/8n$. \pd{In fact, such replacement affects the approximation accuracy in \eqref{OT_primal}, but the corresponding error due to the change of $p,q$ can be made equal to $O(\e)$ with an appropriate choice of $\e'$, and can be made of the same order as the errors coming from the regularization and the rounding algorithm.}

\pd{The regularization parameter $\gamma$ should be sufficiently small to ensure proper approximation in terms of~\eqref{OT_primal}. On the other hand, Theorem~\ref{th:pure_sink} states that the smaller is $\gamma$, the slower is the convergence. The tradeoff is given by the following approximating algorithm from~\cite{dvurechensky2018computational}.}
\begin{algorithm}[H]
    \caption{Approximate OT by Sinkhorn}
    \label{Alg:OTbyS}
    \begin{algorithmic}[1]
        \REQUIRE{Accuracy $\e$.}
        \STATE Set $\e' = \frac{\e}{8\|C\|_\infty}$.
        \STATE 	Set $(\check{p},\check{q})=\left(1 - \frac{\e'}{8}\right)\left((p,q) + \frac{\e'}{8n} (\one_n,\one_n) \right)$
        \STATE Calculate  $\check{\pi}$ s.t. $\|\check{\pi}\one_n - \check p\|_1 + \|\check{\pi}^T\one_n - \check q\|_1 \leqslant \e'/2$ by Algorithm~\ref{Sinkhorn} with marginals $\check{p},\check{q}$ and $\gamma = \frac{\e}{4 \ln n}$
        \STATE Find $\hat{\pi}$ as the projection of $\check{\pi}$ on $U(p,q)$ by Algorithm~\ref{Alg:round}.
        \ENSURE $\hat{\pi}$
  \end{algorithmic}
\end{algorithm}

The followin theorem from~\cite{dvurechensky2018computational} presents a complexity bound for Algorithm~\ref{Alg:OTbyS}.
\begin{theorem}
\label{thm:complexity_OTbyS}
    For a given $\e>0$, Algorithm~\ref{Alg:OTbyS} returns $\hat \pi \in U(p,q)$ s.t.
    \begin{equation*}
         \la C, \hat \pi \ra -  \la C, \pi^* \ra \leqslant \e,
    \end{equation*}
    and requires
    \[
    O\left(\left(\frac{\norm{C}_\infty}{\e}\right)^2 M_{n} \ln n \right)
    \]
    arithmetic operations, where $M_{n}$ is a time complexity of one iteration of Algorithm~\ref{Sinkhorn=AM}.
\end{theorem}

As each iteration of Algorithm~\ref{Sinkhorn=AM} requires a matrix-vector multiplication, the general bound is $M_{n} = O(n^2)$. However, for some specific forms of matrix $C$ it is possible to achieve better complexity, e.g.\ $M_{n} = O(n \log n)$ via FFT \cite{peyre2019computational}.


\subsection{Accelerated Sinkhorn's Algorithm} \label{sec:accsink}

In this subsection, the Primal-Dual AAM algorithm~\cite{guminov2021combination} describes. The algorithm is designed to solve a pair of generally formulated problems~\eqref{primal_general} and~\eqref{dual_general}.
The algorithm solves dual problem~\eqref{dual_general} which should be reformulated as a minimization problem and must have Lipchitz continuous gradient. To obtain this for OT problem one have to suppose $Q =\{ \pi \in \R_+^{n \times n}: \one_n^T \pi\one_n = 1\}$ in~\eqref{primal_general}, since entropic regularizer is $\gamma$ strongly convex on $Q$ w.r.t. $\ell_1$ norm~\cite{nesterov2005smooth}. 

\pd{When applied to OT problems, the idea is to use alternating minimization w.r.t. either $u$ or $v$ in \eqref{OT_dual_lip}.}
Remarkably, Lemma~\ref{Sinkhorn=AM} also holds for the  problem~\eqref{OT_dual_lip}, and partial explicit minimization is possible via the same formulas. 
The authors of~\cite{guminov2021combination} proposed to replace in the classical AGD methods the gradient step with a step of explicit minimization w.r.t. one of the blocks of variables. The resulting algorithm $m$ times theoretically slower than its gradient counterpart, where $m$ is the number of blocks of variables used in alternating minimization. But in practice the algorithm works faster~\cite{guminov2021combination}.

The algorithm performs full-arguments function updates, for that reason it is convenient to introduce a vector containing both dual variables $y$ and $z$. Formally $\mu = \left(y^T, z^T\right)^T$, and notations $\phi(y,z)$ from~\eqref{OT_dual_lip} and $\phi(\mu)$ are equivalent.

\begin{algorithm}[H]
    \caption{Primal-Dual AAM~\cite{guminov2021combination}}
    \label{PDAAM-2}
    \begin{algorithmic}[1]
       \STATE $A_0=a_0=0$, $\eta^0=\null_{2n}$, $\zeta^0:=\mu^0:=\eta^0$, $\check{\pi}^{0} = \pi (\eta^0)$
       \FOR{$k \geqslant 0$}
            \STATE Set $\beta_k = \argmin\limits_{\beta\in [0,1]} \phi\left(\eta^t + \beta (\zeta^t - \eta^t)\right)$
            \STATE Set $\mu^{t}=\beta_k \zeta^t+(1-\beta_k)\eta^t$
            \STATE Choose $\xi=\argmax\limits_{\xi\in\{u, v\}} \|\nabla_\xi \phi(\mu^t)\|_2^2$
            \STATE Set $\eta^{t+1}=\argmin\limits_{\xi} \phi(\xi)$
            \STATE Find $a_{k+1}$, $A_{k+1} = A_{k} + a_{k+1}$  from  \[\phi(\mu^t)-\frac{a_{k+1}^2}{2(A_k+a_{k+1})}\|\nabla \phi(\mu^t)\|_2^2=\phi(\eta^{t+1})\]\\
            \STATE Set $\zeta^{t+1} = \zeta^{t} - a_{k+1}\nabla \phi(\mu^t)$
            \STATE Set $\check{\pi}^{t+1} = \frac{a_{k+1}\pi(\mu^{t})+A_k\check{\pi}^{t}}{A_{k+1}}$
        \ENDFOR
    \ENSURE The points $\check{\pi}^{t+1}$, $\eta^{t+1}$
    \end{algorithmic}
\end{algorithm}

\begin{theorem} \label{PD-bounds}
    Let the objective $g(\pi)$ in problem~\eqref{primal_general} be $\gamma$-strongly convex w.r.t. $\|\cdot\|_E$, and let $\|\mu^*\| \leqslant D$. Then, for the sequences $\check{\pi}^{t},\eta^{t}$, $t\geqslant 0$, generated by Algorithm \ref{PDAAM-2}, 
    \begin{align}
        &|\phi(\eta^t) + g(\check{\pi}^t)| \leqslant\frac{8m\|\bm{A}\|^2_{E\to H}D^2}{\gamma t^2}, 
        \\
        &\|\bm{A} \check{\pi}^t - b \|_H \leqslant \frac{8m\|\bm{A}\|^2_{E\to H}D}{\gamma t^2},
        \\
        &\|\check{\pi}^k-\pi_{\gamma}^*\|_E \leqslant \frac{4m  \|\bm{A}\|_{E\to H}D}{\gamma t},
        \label{eq:APDAAM_bound}
    \end{align}
    where $\|\bm{A}\|_{E\to H}$ is the norm of $\bm{A}$ as a linear operator from $E$ to $H$, i.e. $\|\bm{A}\|_{E \rightarrow H} = \max_{x,y} \left\{{\la Ax, y\ra : \|x\|_E = 1, \|y\|_H = 1}\right\}$, and $\|\cdot\|_H = \|\cdot\|_2$.
\end{theorem}

To apply Algorithm~\ref{PDAAM-2} for solving problem~\eqref{def:regoptimaltransport} the authors of~\cite{guminov2021combination} equipped the primal space $E$ with 1-norm and the dual space $H$ with 2-norm. This for the case of OT problem implies that $\|\bm{A}\|^2_{E\to H} = \|\bm{A}\|^2_{1\to2}=2$. They applied the algorithm for minimizing the dual function written as in~\eqref{OT_dual_lip}. Particularly, the authors of~\cite{guminov2021combination}  suppose further that change of variables~\eqref{var_change}  is done, and the dual objective can be minimized explicitly in $u$ and $v$ as mentioned above.
The notation for partial gradient of~\eqref{OT_dual_lip} is introduced as follows
\[
    \nabla_u\varphi(\lambda) = \partial \varphi(u,v) /\partial u = \gamma \left( p-\frac{B(u,v)\one_n}{\one_n^TB(u,v)\one_n} \right)
\]
\[
\nabla_v\varphi(\lambda) = \partial \varphi(u,v) /\partial v = \gamma \left( q-\frac{B(u,v)^T\one_n}{\one_n^TB(u,v)\one_n} \right).
\]


Theoretical analysis of the algorithm relies on the choice of the starting point $\eta_0=\null_{2n}$, that allows to estimate the distance to the solution $\mu^*$ of~\eqref{OT_dual_lip} as
\[
    \|\mu^*\|_2\leqslant D := \sqrt{n/2}\left(\|C\|_\infty-\frac{\gamma}{2}\ln\left(\min\left\{\min_i p_i, \min_j q_j\right\}\right)\right).
\]

Computations of the objective and its gradient should  also be implemented using numerically stable computation of the log of the sum of exponents.

The authors of~\cite{gasnikov2015universal} pointed out that the objective~\eqref{OT_dual_lip}, its gradient and equation~\eqref{eq:pvr} are invariant under transformations $u\to u+t_u\one_n$, $v\to v+t_v\one_n$, with $t_u,t_v\in\mathbb{R}$. In practice such variable transformation (with $t_\xi = -\|\xi\|_\infty$, $\xi\in \{u,v\}$) performed after steps 4 and 8 of Algorithm~\ref{PDAAM-2} can lead to better numerical stability when $\gamma$ is small.

As for Sinkhorn's algorithm, parameter $\gamma$ should be sufficiently small to ensure proper approximation in terms of~\eqref{OT_primal}. On the other hand, Theorem~\ref{PD-bounds} states that the smaller is $\gamma$, the slower is the convergence. The tradeoff is given by the Algorithm~\ref{Alg:OTbyAccS} from~\cite{guminov2021combination}.

\begin{algorithm}[ht!]
  \caption{Accelerated Sinkhorn for OT}
  \label{Alg:OTbyAccS}
  \begin{algorithmic}[1]
    \REQUIRE{Accuracy $\e$.}
    \STATE Set $\gamma = \frac{\e}{3 \ln n}$, $\e' = \frac{\e}{8 \|C\|_\infty}$.
    
    \STATE 	Set $(\check{p},\check{q})=\left(1 - \frac{\e'}{8}\right)\left((p,q) + \frac{\e'}{8n} (\one,\one) \right)$
    \FOR{$k=1,2,...$}
        \STATE Perform an iteration of Algorithm~\ref{PDAAM-2} for the OT problem with marginals $\check{p},\check{q}$ and calculate $\check{\pi}^t$ and $\eta^t$.
        \STATE Find $\hat{\pi}$ as the projection of $\check{\pi}^t$ on $U(p,q)$ by Algorithm~\ref{Alg:round}.
        \STATE {\textbf{if} $\la C,\hat{\pi}-\check{\pi}^t\ra \leqslant \frac{\e}{6}$ and $g(\check{\pi}^t)+\varphi(\eta^t) \leqslant \frac{\e}{6}$}
        \STATE \textbf{then} Return $\hat{\pi}$
   \ENDFOR
  \end{algorithmic}
\end{algorithm}


Primal-dual property of Algorithm~\ref{PDAAM-2} allows one to choose greater values of regularization parameter $\gamma$ than in Algorithm~\ref{Alg:OTbyS}. The latter also leads to numerical stable behavior of the algorithm.

Next theorem from~\cite{guminov2021combination} presents the complexity bound for Algorithm~\ref{Alg:OTbyAccS}.
\begin{theorem}
\label{thm:complexity_OTbyAccS}
    For a given $\e>0$, Algorithm~\ref{Alg:OTbyAccS} returns $\hat \pi \in U(p,q)$ s.t.
    \begin{equation*}
         \la C, \hat \pi \ra -  \la C, \pi^* \ra \leqslant \e,
    \end{equation*}
    and requires
    \[
        O\left(n^{5/2}\sqrt{\ln n}\frac{\norm{C}_\infty}{\e} \right)
    \]
    arithmetic operations.
\end{theorem}

\section{Wasserstein Barycenters}
One of the most fundamental questions in statistics and data analysis is inferring
the mean of a random variable given its realizations $p_1, \cdots, p_m$.
Although this is conceptually straightforward when the data lie in a Euclidean space or a Hilbert space, many datasets exhibit complex geometries that are far from Euclidean~\cite{heinemann2020randomised}.
A Fr{\'e}chet mean (or barycenter)~\cite{frechet1948elements}  extends
 the Euclidean mean $q = \argmin_q \sum_l\|p_l - q\|^2$ of points   $p_1, \cdots, p_m \in S_n(1)$   to a non-Euclidean case by replacing the norm $\|p_l - q\|$ with an arbitrary distance function. 
Using this notion and optimal transport distances,  we define fixed-support Wasserstein barycenter (WB) of probability measures $p_1, \dots, p_m$ as a solution to the following convex optimization problem:
\begin{equation}\label{prob:unreg_bary}
    \min_{q \in S_n(1)} \frac{1}{m} \sum_{l=1}^m \W(p_l, q).
\end{equation}
The existence and consistency of WB were studied in \cite{le2017existence}. One of the major application  is  deformation models: an unknown template distribution is warped from different observations by a random deformation process. Then the barycenter of the observations is a proper estimate for the unobserved template (it is a consistent estimate under some conditions on deformation process).
Moreover, WB are used in Bayesian computations \cite{srivastava2015wasp}, texture mixing \cite{rabin2011wasserstein}, clustering  ($k$-means for probability measures) \cite{del2019robust}, shape interpolation and color transferring \cite{Solomon2015} and neuroimaging \cite{gramfort2015fast}.  

Generally, problem \eqref{prob:unreg_bary} cannot be solved analytically, however, for instance, in case of 
Gaussian measures, it has an explicit solution and the barycenter  is also Gaussian. For arbitrary probability measures, problem \eqref{prob:unreg_bary} is computationally expensive. The entropic regularization approach was extended for barycenters and leads to the following problem formulation
\begin{equation}\label{prob:reg_bary}
    \min_{q \in S_n(1)} \frac{1}{m} \sum_{l=1}^m \W_\gamma(p_l, q).
\end{equation} 
Similarly, to the entropy-regularized OT problem, which can be solved by the Sinkhorn’s algorithm, for problem  \eqref{prob:reg_bary}the iterative Bregman projections (IBP) algorithm  \cite{benamou2015iterative} was designed. The IBP  is an extension of the Sinkhorn’s algorithm for $m$ measures,  and, hence, its complexity is  $m$ times larger than the Sinkhorn's algorithm complexity. Namely, it is $ \widetilde O\left({ mn^2  \|C\|^2_\infty}/{(\gamma \e)}\right)$ \cite{kroshnin2019complexity}. 
An analog of the accelerated Sinkhorn's algorithm for the WB problem of $m$ measures is the  accelerated IBP algorithm with complexity $\widetilde O\left({mn^2 \sqrt{n} \|C\|_\infty}/{\sqrt{\gamma \e}}\right)$  \cite{guminov2019accelerated}, that is also $m$ times larger than the accelerated Sinkhorn complexity. 
Another fast version of the IBP algorithm was  proposed in \cite{lin2020fixed}, named FastIBP algorithm  with the complexity $ \widetilde O\left({mn^2\sqrt[3]{n}  \|C\|^{4/3}_\infty}/{(\gamma \e)^{2/3}}\right)$. Problem \eqref{prob:reg_bary} can be also solved  by primal-dual algorithm (primal-dual accelerated gradient descent) studied in \cite{dvurechenskii2018decentralize,dvinskikh2019primal}   with  the complexity $\widetilde O \left(mn^{5/2}\|C\|_\infty /\sqrt{\gamma \e}\right)$. Primal-dual first-order methods rely on the fact  that there exists closed-form representations for the dual function for $\W_\gamma (p,q)$ defined by the Fenchel–Legendre transform  \cite{agueh2011barycenters,cuturi2016smoothed}:
\begin{align*}
     \W_{\gamma, p}^*(u) 
     &=  \max_{ q \in S_1(n)}\left\{ \la u, q \ra - \W_{\gamma}(q, p) \right\} \notag \\
  &= 
    \gamma\sum_{j=1}^n [p]_j \log\left( \sum_{i=1}^n \exp\left(\frac{[u]_i - C_{ji}}{\gamma}\right) \right)-\gamma\la p,\log p\ra,
\end{align*}
  where   $[p]_j$ and $[u]_i$ are the $j$-th and $i$-th components   of  $p$ and $u$ respectively. 
  The gradient of  $\W_{\gamma, p}^*(u) $ is Lipschitz continuous and has also a closed-form solution \cite[Theorem 2.4]{cuturi2016smoothed}
 \begin{align}\label{eq:primal_sol_recov3323}
 [\nabla \W^*_{\gamma,p} (u)]_l = \sum_{j=1}^n [p]_j \frac{\exp\left(([u]_l-C_{lj})/\gamma\right)  }{\sum_{\ell=1}^n\exp\left(([u]_\ell-C_{\ell j})/\gamma\right)},
\end{align}
for all $ l =1,...,n$.



However,  all  entropy-regularization-based approaches are numerically unstable when the regularization parameter $\gamma$ is  small (this also means that the precision $\e$ is small as $\gamma $ must be selected proportionally to $\e$ \cite{peyre2019computational,kroshnin2019complexity}). The authors of \cite{dvinskikh2021improved} solve  problem \eqref{prob:unreg_bary}  without additional regularization by rewriting it  as a saddle-point problem. 
Their saddle-point approach for the WB problem relies on the fact that non-regularized optimal transport problem \eqref{def:optimaltransport} has  a bilinear  saddle-point representation
 \cite{jambulapati2019direct}:
\begin{equation*}
    \resizebox{\linewidth}{!}{
    $\W(p,q)=\min_{x \in S_{n^2}(1)} \max_{y\in [-1,1]^{2n}} \left\{ \la d, x \ra +2\|d\|_\infty\left(~ y^\top Ax - \left\la \begin{pmatrix}
        p\\
        q
    \end{pmatrix}, y \right\ra\right)\right\}$.
    }
\end{equation*}

Here $d$ is the vectorized cost matrix $C$, $x \in S_{n^2}(1)$ is the vectorized transport plan $\pi$, and \[A =  \begin{pmatrix}
    I_{n\times n} &\otimes &\boldsymbol{1}^\top_{n}\\
     \boldsymbol{1}^\top_{n} &\otimes &I_{n\times n}
 \end{pmatrix} = \{0,1\}^{2n\times n^2} \] is the incidence matrix, where $\otimes$ is the Kronecker product. 
The proposed algorithm \cite{dvinskikh2021improved}  has a near-optimal rate $\widetilde O \left(mn^{2}\|C\|_\infty /\e\right)$ and it is an extension of the saddle-point algorithm from \cite{jambulapati2019direct}  for saddle-point formulation of optimal transport problem.   

\subsection{Approximating Non-regularized WB by IBP}
In this subsection, we consider the Iterative Bregman Projections algorithm~\cite{benamou2015iterative} to solve the problem~\eqref{prob:reg_bary}. We slightly reformulate this problem as
\begin{multline}\label{prob:reg_bary_2}
    \min_{q \in S_n(1)} \frac{1}{m}\sum_{l=1}^m \W_\gamma (p_l,q) = \min_{\substack{q \in S_n(1), \\ \pi_l \in U(p_l, q),\\
    l = 1, \dots, m }} \frac{1}{m} \sum_{l = 1}^m  \bigl\{\la \pi_l, C \ra + \gamma H(\pi_l) \bigr\} 
    \\ 
    = \min_{\substack{ \pi_l \in \R_+^{n \times n}\\ \pi_l \one = p_l,\; \pi_l^T \one = \pi_{l+1}^T \one, \\
    l=1,\dots,m}} \frac{1}{m}\sum_{l = 1}^m   \bigl\{\la \pi_l, C \ra + \gamma H(\pi_l)  \bigr\}.
\end{multline}
The dual problem to~\eqref{prob:reg_bary_2} is (up to a multiplicative constant)
\begin{equation} \label{prob:dual_2}
    \min_{\substack{u, v\\ \sum_l v_l = 0}}  f(u,v) := \frac{1}{m} \sum_{l = 1}^m \bigl\{\la \one_n, B(u_l, v_l) \one_n \ra - \la u_l, p_l \ra\bigr\},
\end{equation}
$u = [u_1, \dots, u_m]$, $v = [v_1, \dots, v_m]$, $u_l, v_l \in \R^n$, and
\begin{equation}\label{eq:B_def}
    B(u_l, v_l) := \diag\left(e^{u_l}\right) K \diag\left(e^{v_l}\right), \quad
    K := \exp\left(-\frac{C}{\gamma}\right).
\end{equation}
The primal variable can be reconstructed by the formula $\pi_l(u_l,v_l) = B(u_l,v_l)$.

The following lemma shows that the dual problem~\eqref{prob:dual_2} can also be minimized exactly over the variables $u,v$.

\begin{lemma} \label{lemma:IBP_explicit_step}
    The iterations of IBP algorithm
    \[u^{t+1} \in \argmin_{u\in\mathbb{R}^{m\times n}}\vp(u,v^t),\quad v^{t+1}\in \argmin_{v\in\mathbb{R}^{m\times n}, \,\, \sum_l v_l = 0}\vp(u^{t+1},v)\] can be written explicitly as
    \begin{equation} \label{eq:u_cond}
        u^{t+1}_l 
        = u^t_l + \ln p_l - \ln \left(B\left(u_l^t, v_l^t\right) \one_n\right) 
        = \ln p_l - \ln K e^{v_l^t},
    \end{equation}
    and 
    \begin{equation} \label{eq:v_cond}
        v^{t + 1}_l
        = v^t_l + \ln q^{t + 1} - \ln q^t_l 
        = \frac{1}{m} \sum_{k = 1}^m \ln K^T e^{u_k^t} - \ln K^T e^{u_l^t},
    \end{equation}
    where $q^t_l := B^T(u^t_l, v^t_l) \one_n$, $q^{t + 1} := \exp\left(\frac{1}{m}\sum_{l = 1}^m \ln q^t_l\right)$.
\end{lemma}

To solve the dual problem, the authors of~\cite{kroshnin2019complexity} reformulate the IBP algorithm as Algorithm~\ref{Alg:dual_IBP} allowing to solve simultaneously the primal and dual problem and have an adaptive stopping criterion (see line 7), which does not require to calculate any objective values.   

\begin{algorithm}[H]
    \caption{Iterative Bregman Projections (IBP) }
    \label{Alg:dual_IBP}
    \begin{algorithmic}[1]
        \REQUIRE $C$, $p_1, \dots, p_m$, $\gamma > 0$, $\e' > 0$
        \STATE $u_l^0 := 0$, $v_l^0 := 0$, $l = 1, \dots, m$ 
        \REPEAT {}
            \STATE   $v_l^{t + 1} := \frac{1}{m} \sum_{k = 1}^m \ln K^T e^{u_k^t} - \ln K^T e^{u_l^t}, \quad u^{t + 1} := u^t$ 
            \STATE    $t := t + 1$
            \STATE  $u^{t + 1}_l := \ln p_l - \ln K e^{v_l^t}, \quad v^{t + 1} := v^t$ 
            \STATE $t := t + 1$
        \UNTIL{$\frac{1}{m} \sum_{l = 1}^m  \norm{B^T(u_l^t, v_l^t) \one_n - \bar{q}^t}_1 \le \e'$, \\ \quad where $\bar{q}^t := \frac{1}{m} \sum_{l = 1}^m  B^T(u_l^t, v_l^t) \one_n$}
        \ENSURE $B(u^t_1, v^t_1), \dots, B(u^t_m, v^t_m)$  
    \end{algorithmic}
\end{algorithm}

To find an approximate solution to the initial problem~\eqref{prob:unreg_bary} the authors of~\cite{kroshnin2019complexity} applied Algorithm~\ref{Alg:dual_IBP} with a suitable choice of $\gamma$ and $\e'$. For simplicity we average marginals $q_1, \dots, q_m$ with uniform weights, which leads to Algorithm~\ref{Alg:IBP_to_bary}.

\begin{algorithm}[H]
    \caption{Finding Wasserstein barycenter by IBP~\cite{kroshnin2019complexity}}
    \label{Alg:IBP_to_bary}
    \begin{algorithmic}[1]
        \REQUIRE {Accuracy $\e$; cost matrices $C_1, \dots, C_m$; marginals $p_1, \dots, p_m$}
        \STATE {Set  $\gamma := \frac{1}{4} \frac{\e}{\ln n}, \quad \e' :=  \frac{1}{4} \frac{\e}{\norm{C}_\infty}$}
        \STATE {Find $\check{B}_l := B(u^t_l, v^t_l)$, $l = \{1, \dots, m\}$ by Algorithm~\ref{Alg:dual_IBP} with accuracy $\e'$}
        \STATE {$\bar q := \sum_{l = 1}^m  \check{B}_l^T \one_n / \sum_{l = 1}^m \la \one_n, \check{B}_l \one_n \ra$}
        \STATE Find $\hat{\pi_l}$ as the projection of $\check{B}_l$ on $U(p_l,\bar q)$ by Algorithm~\ref{Alg:round} for $l = \{1, \dots, m\}$.
        \ENSURE {$\bar q, \hat{\pi}_1, \dots, \hat{\pi}_m$}
    \end{algorithmic}
\end{algorithm}

The following theorem presents the complexity bound obtained in~\cite{kroshnin2019complexity} for Algorithm~\ref{Alg:IBP_to_bary}.
\begin{theorem}
\label{thm:complexity_IBP_to_bary}
    For a given $\e>0$, Algorithm~\ref{Alg:IBP_to_bary} returns $q \in S_n(1)$ s.t.
    \[
    \frac{1}{m} \sum_{l = 1}^m  \W(p_l, \bar q) 
    - \frac{1}{m} \sum_{l = 1}^m \W(p_l, q^*) \le \e,
    \]
    where $q^*$ is a solution to a non-regularized problem~\eqref{prob:unreg_bary}.
    It requires
    \[
    O\left(\left(\frac{\norm{C}_\infty}{\e}\right)^2 M_{m,n} \ln n + m n\right)
    \]
    arithmetic operations, where $M_{m,n}$ is a time complexity of one iteration of Algorithm~\ref{Alg:dual_IBP}.
\end{theorem}

As each iteration of Algorithm~\ref{Alg:dual_IBP} requires $m$ matrix-vector multiplications, the general bound is $M_{m,n} = O(m n^2)$. However, for some specific forms of matrix $C$ it is possible to achieve better complexity, e.g.\ $M_{m,n} = O(m n \log n)$ via FFT \cite{peyre2019computational}.

\subsection{Accelerated IBP Algorithm}
Additional tautological constraint $\one_n^T \pi_l \one_n = 1$, $l=1,\dots,m$  were introduced in~\cite{guminov2021combination} in the problem~\eqref{prob:reg_bary_2}, leading to the dual problem with the objective function having Lipchitz continuous gradient:

\begin{equation}
    \label{WB_dual_supp}
       \min_{
       \substack{
        u,v\\
        \sum_{l=1}^m v_l=0
       }
       }
       \varphi(u, v) := \frac{\gamma}{m}\sum\limits_{l=1}^{m} \left\{ \ln\left(\mathbf{1}^T B\left(u_{l}, v_{l}\right)\mathbf{1}\right)-\left\langle u_{l}, p_{l}\right\rangle\right\},
\end{equation}
where $u = [u_1,\ldots , u_m], v = [v_1, \ldots , v_m], u_l,v_l\in\mathbb{R}^n$. For this problem the primal-dual relation between the variables is given by
\begin{equation}  
    \pi_l(u_l,v_l) = B(u_l,v_l) / \one_n^T B(u_l,v_l)\one_n.
\end{equation}

Lemma~\ref{lemma:IBP_explicit_step} holds also for the  problem~\eqref{WB_dual_supp}, and  explicit alternating minimization is possible. This allows one to apply Algorithm \ref{PDAAM-2} to approximate the solution of WB problem with accuracy $\e$. Algorithm~\ref{Alg:WBbyAccS} is the pseudocode of such algorithm proposed in~\cite{guminov2021combination}.

\begin{algorithm}[H]
  \caption{Accelerated IBP~\cite{guminov2021combination}}
  \label{Alg:WBbyAccS}
  \begin{algorithmic}[1]
    \REQUIRE{Accuracy $\e$.}
    \STATE Set $\gamma = \frac{\e}{2 \ln n}$, $\e' = \frac{\e}{8\|C\|_\infty}$.
    
    \STATE 	Set $\check{p_l}=\left(1 - \frac{\e'}{4}\right)\left(p_l + \frac{\e'}{4n} \one \right)$
    \FOR{$k=1,2,...$}
    
        \STATE Perform an iteration of Algorithm~\ref{PDAAM-2} for the WB problem with marginals $\check{p}_l$ and calculate $\check{\pi}_l^k$, $l=1,\dots,m$ and $\eta^k$
        \STATE Find $\bar q = \frac{1}{m}\sum_{l = 1}^m (\check{\pi}_l^k)^T \one$
        \STATE Calculate $\hat{\pi}_l$ as the projection of $\check{\pi}_l^k$ on $U(\check{p}_l, \bar q)$ by Algorithm~\ref{Alg:round}
        \STATE {\textbf{if}  $\frac{1}{m} \sum_{l = 1}^m \left\{\la C,\hat{\pi}_l\ra - \la C,\check{\pi}_l^k \ra \right\} \leqslant  \frac{\e}{4}$ and $\frac{1}{m}\sum_{l = 1}^m   \bigl\{\la \check{\pi}_l, C \ra + \gamma H(\check{\pi}_l)  \bigr\}+\varphi(\eta_k) \leqslant \frac{\e}{4}$}
        \STATE \textbf{then} Return $\bar q, \hat{\pi}_1, \dots, \hat{\pi}_m$.
   \ENDFOR
  \end{algorithmic}
\end{algorithm}


Since each iteration requires $O(mn^2)$ arithmetic operations, which is the same as in the IBP algorithm, the authors of~\cite{guminov2021combination} obtained the following theorem that gives the complexity bound for Algorithm~\ref{Alg:WBbyAccS}:
\begin{theorem}
\label{thm:wbbyaccs}
    For a given $\e>$, Algorithm~\ref{Alg:WBbyAccS} returns $\bar q \in S_n(1)$ s.t.
    \[
    \frac{1}{m} \sum_{l = 1}^m  \W(p_l, \bar q) 
    - \frac{1}{m} \sum_{l = 1}^m \W(p_l, q^*) \le \e,
    \]
    where $q^*$ is a solution to a non-regularized problem~\eqref{prob:unreg_bary}.
    It requires
    \[O\left(\frac{mn^{5/2}\sqrt{\ln n}\|C\|_\infty}{\varepsilon}\right).\]
    arithmetic operations.
\end{theorem}

Numerical stability features described in Sections~\ref{sec:sink} and~\ref{sec:accsink} should be applied to IBP and accelerated IBP algorithms. Note that~\cite{kroshnin2019complexity} did not separate $p_l$ entries from zero, by replacing input measures $p_l$ with $\check{p}_l$. The latter can lead to numerical instability of the IBP algorithm. But of course it's possible to prove similar complexity bounds with the replacement.  Moreover, accelerated IBP allows the use of a greater value of regularization parameter $\gamma$, which leads to better numerical stability.

\subsection{Distributed  Computation of Barycenters}

Some practical  applications involve large quantities of
data without centralized storage: data  can be stored or collected in a distributed manner (e.g., sensors in a sensor network obtain the state of the environment from different geographical parts, or  micro-satellites collecting local information).
Thus, a question of distributed or even decentralized computation of the  WB arises. To this end, we consider the  WB problem over a computational network, i.e., a system of
computing nodes (agents, machines, processing units), whose interactions are constrained. They process faster and  more data than one computer since the work is divided  between many computing nodes.  Each agent privately holds a  probability distribution, can sample from it, and seeks to cooperatively compute the barycenter of all probability measures exchanging the information with its neighbors.

Decentralized formulations of the WB problem 
both for the saddle-point and dual approaches
are based on introducing at each node $i$ a local copy $q_i$ of the global variable and 
introducing an artificial constraint
$q_1= ...=q_m \in \R^n$ which is further replaced with affine constraint  $  \boldsymbol W  \boldsymbol q=0$ (in the saddle-point approach)  and $\sqrt{\boldsymbol W}   \boldsymbol q = 0$ (in the dual approach), where 
$\boldsymbol q = (q_1^\top, ..., q_m^\top)^\top$ is column vector and 
$\boldsymbol W$ is referred as the communication matrix of a decentralized system:
$\boldsymbol W = W \otimes I_n$, where  $I_n$ is the identity matrix and $W{\in \mathbb{R}^{m\times m}}$ is the   Laplacian matrix
 of the graph $\mathcal{G}= (V, E)$ representing the computational network:
\begin{align*}
[W]_{ij} = \begin{cases}
-1,  & \text{if } (i,j) \in E,\\
\text{deg}(i), &\text{if } i= j, \\
0,  & \text{otherwise.}
\end{cases}
\end{align*}
Here $\text{deg}(i)$ is the degree of the node $i$, i.e., the number of neighbors of the node. 
From the definition of matrix $ \boldsymbol W$ it follows that 
\[\sqrt{\boldsymbol W}  \boldsymbol q=0 \Longleftrightarrow  \boldsymbol W \boldsymbol q = 0 
\Longleftrightarrow q_1 = ... = q_m.\] 
The affine constraint $ \boldsymbol W \boldsymbol q=0$ (or $\sqrt{\boldsymbol W }\boldsymbol q = 0$) is brought to the objective via the Fenchel--Legendre transform.
Thus,  to state the regularized Wasserstein barycenter problem \eqref{prob:reg_bary} in a decentralized manner, we rewrite it as follows
\begin{equation}\label{eq:W_bary_reg22}
   \min_{\substack{q_1=...=q_m, \\q_1,...,q_m \in S_n(1)} } \frac{1}{m} \sum_{i=1}^m {\W}_\gamma(p_i,q_i) = \min_{\substack{\sqrt{\boldsymbol W} \boldsymbol q=0, \\q_1,...,q_m \in S_n(1)} } \frac{1}{m} \sum_{i=1}^m {\W}_\gamma(p_i,q_i).
\end{equation}
The dual problem to \eqref{eq:W_bary_reg22} is
\begin{equation}\label{eq:W_bary_regdual}
  \min_{\boldsymbol y \in \R^{nm}} \left\{\W^*_{\gamma,\boldsymbol p}\left(\sqrt{\boldsymbol W}\y\right)= \frac{1}{m}\sum_{i=1}^{m} {\W}^*_{\gamma, p_i}\left(m\left[\sqrt{\boldsymbol W}\boldsymbol y\right]_i\right)\right\},
\end{equation}
where ${\W}^*_{\gamma, p_i}\left(m\left[\sqrt{\boldsymbol W}\boldsymbol y\right]_i\right)$ is the Fenchel-Legendre transform, $\p = (p_1^\top, \cdots, p_m^\top)^\top$, and $\y = (y_1^\top, \cdots, y_m^\top)^\top \in \R^{nm}$ is the Lagrangian dual multiplier.

As the primal function  ${\W}_\gamma$ is strongly convex in the $\ell_2$-norm on $S_n(1)$, then the dual function ${\W}^*_{\gamma, \p}$ is  $L$-Lipschitz smooth, or has Lipschitz continuous gradient. The constant $L$ for ${\W}^*_{\gamma, \p}(\sqrt{\boldsymbol W}\y)$ is defined via communication matrix $\boldsymbol W$ and regularization  parameter $\gamma$. Hence, accelerated gradient descent-based method  can be used, which is optimal in terms of the number of iterations and oracle calls. 
For simplicity, a decentralized procedure  solving dual problem \eqref{eq:W_bary_regdual} can be demonstrated on  the gradient descent as follows
\begin{equation*}
\y^{k+1}= \y^{k}- \frac{1}{L} \nabla{\W}^*_{\gamma, \q}(\sqrt{\boldsymbol W}\y^k) = \y^{k} -\frac{1}{L}\sqrt{\boldsymbol W}\q(\sqrt{\boldsymbol W}\y^k). \end{equation*}
Without change of variable, it is unclear how to execute this  procedure in a distributed fashion. Let  $\u := \sqrt{\boldsymbol W}\y $, then the gradient step  multiplied by $\sqrt {\boldsymbol W}$ can be rewritten as 
\[\u^{k+1} = \u^{k} -\frac{1}{L}\boldsymbol W\q(\u^{k}),\]
 where $[\q(\boldsymbol u)]_i = q_i(u_i) = \nabla \W^*_{\gamma,p_i}(u_i)$ from \eqref{eq:primal_sol_recov3323}, $i=1,...,m$.
This procedure can be performed in a decentralized manner on a distributed network.
The vector $\boldsymbol W\q(\boldsymbol u)$ naturally defines communications with neighboring nodes due to the structure of communication matrix $\boldsymbol W$ as  the elements of communication matrix are zero for non-neighboring nodes. 
Moreover, 
in the dual approach which is based on gradient method, the randomization of $\nabla \W_{\gamma, p_i}^*(u_i)$ can be used  to reduce the complexity of calculating the true gradient, that is $O(n^2)$ arithmetic operations, by calculating its stochastic approximation of $O(n)$ arithmetic operations. 
The randomization for the true gradient \eqref{eq:primal_sol_recov3323} is achieved   by taking the $j$-th term in the sum with probability $[p]_j$ 
\[
[\nabla \W_{\gamma,p}^*(u,\xi)]_l =  \frac{\exp\left(([u]_l-C_{l \xi})/\gamma\right)  }{\sum_{\ell=1}^n\exp\left(([u]_\ell-C_{\ell\xi })/\gamma\right)}, \,\, \forall l =1,...,n.
\]
where we replaced index $j$ by $\xi$ to underline its randomness.
This is the motivation for considering the first-order methods with stochastic oracle. 

Algorithms solving the regularized  WB problem \eqref{prob:reg_bary} such as the IBP \cite{benamou2015iterative}, accelerated IBP \cite{guminov2019accelerated}   and the FastIBP \cite{lin2020fixed}  can be implemented in a distributed manner but with a central fusion center that coordinates the actions of the parallel machines. 
Primal-dual algorithm based on the accelerated gradient descent \cite{uribe2017optimal,dvurechenskii2018decentralize,dvinskikh2019primal}  and decentralized mirror prox (DMP) algorithm \cite{rogozin2021decentralized} can be implemented in a decentralized manner on an arbitrary connected network. 
The oracle complexities of all these aforementioned (non-randomized) methods are presented in Table \ref{Tab:distr_comp}.
Communication complexities are obtained by multiplying the oracle complexities by the square root of the condition number $\chi = \frac{\lm_{\max} (W)}{\lm^+_{\min} (W)}$,  where $\lm_{\max}(W)$ and $\lm^+_{\min}(W)$ are the maximum and minimum
non-zero  eigenvalues of matrix $W$ respectively.
\dm{\begin{table}[ht!]
    \centering
    \caption{ Distributed algorithms  for the WB problem}
    \resizebox{\columnwidth}{!}{%
    \begin{threeparttable}
    \begin{tabular}{ll}
    Algorithm   &     Complexity per node$^{(1)}$   \\
     \midrule
     IBP \cite{benamou2015iterative,kroshnin2019complexity}    &  
    $n^2/(\gamma\e) $    \\
    Accelerated IBP  \cite{guminov2019accelerated}               & 
       $ n^2 \sqrt{n}/\sqrt{\gamma \e} $ 
       \\
      FastIBP    \cite{lin2020fixed}       & 
    $n^2\sqrt[3]{n} /(\gamma \e )^{2/3}$   
    \\
    Primal-dual algorithm \cite{uribe2017optimal,dvurechenskii2018decentralize,dvinskikh2019primal} &  $ n^2\sqrt{ n }/\sqrt{\gamma\e}$ \\
     {DMP} \cite{rogozin2021decentralized}    & 
       $ n^2\sqrt{ n}  / \e $
       \\
    \end{tabular}
    \begin{tablenotes}
    \item[]  $^{(1)}$ The bounds are obtained by using the Chebyshev acceleration, see e.g. \cite{gorbunov2022recent} and references there in.
    \end{tablenotes}    
    \end{threeparttable}}
    \label{Tab:distr_comp}
\end{table}}

\section*{Acknowledgment}
This research was funded by Russian Science Foundation (project 18-71-10108).

\bibliographystyle{plain}
\bibliography{references}

\end{document}